\documentclass[letterpaper, 12pt]{article}
\title{Set Theory Axioms Using ``Explication''}
\author{E. R. Akemann\\
                  At Home}
\date {August 22, 2017}
\begin {document}
\maketitle
 \section{Overview}
\setlength{\parskip} {1ex} 
In GB (System of G\"{o}del-Bernays), one starts with ``class''  variables (see Appendix A) and a unary relation M(X) meaning X is a set.  The basic axiom for sethood is $ Y \in X \rightarrow M(Y)$.  We introduce a second unary relation SSC(X) meaning that X is a ``set-sized'' class, thus $Y \in X \rightarrow SSC(Y)$. These ``set-sized'' classes will then be subjected to a filtering process, which will determine which SSC's are sets, i.e. explicated $(E(x))$, and which are not. As in GB, a class variable in lower case is automatically and implicitly an SSC.
A dot above $\dot{x}$ denotes an SSC not yet explicated; if already a set, the dot disappears.  Thus, $x \equiv \dot{x} \ \& E(x)$.  It will be clear in the discussion of the axioms, which SSCs become sets and which are filtered out.

\noindent The basic idea of this axiom scheme (called EX) is that the SSCs which become sets are explicated by a \textbf{finite} sequence of ``explicating axioms'', which closely follow the Axioms of ZF (System of Zermelo-Fraenkel), restricted to explicated SSCs. The goal of the filtering process is to filter \textbf{out} any SSC which doesn't explicitly say, in a finite statement, what sets are its members.  The mixing of GB and ZF is necessary. EX starts with GB, because ZF has only sets, while the SSC's are needed to to make EX comprehensible. ZF is needed because the function making axioms of GB allow a non-explicated SSC to become explicated, while ZF's can be structured so as not to allow that to happen.

\noindent See reference COHEN[1] for a complete treatment of GB and ZF.
\newline The following suite of axioms is a \textbf{complete} treatment of EX.  Items which may or may not be true in other axiom systems are not significant.
\section{Axiom of Extension}
$ x = y \iff \forall u (u \in x \iff u \in y)$.\ \ The Axiom of Extension reinforces the overall idea of EX in that an SSC which can't explicitly say, in a finite statement, what elements are members, can't be compared to another SSC, and is thus a poor candidate for being a set.
\section{Set Explicating Axioms}
\begin{enumerate}
\item Null Set \newline 
$\exists \dot{x} (\forall y (y \notin \dot{x} ) ) \ \& \  E(x, NS)\ [x \ denoted \ by \ \phi] $ 
\item Unordered Pairs \newline 
$\forall x,y \  \exists \dot{z} \ \forall w (w \in \dot{z} \iff (w = x) \lor (w=y)) \rightarrow E(z,UP,x,y)$  
\item Sum Set \newline
$\forall x \ \exists \dot{y} ( \forall z(z \in \dot{y} \iff  \exists t(z \in t \ \& \ t  \in x))) \rightarrow E(y,SS,x) $
\item Infinity \newline
$\exists \dot{x} (\phi \in \dot{x} \ \& \  \forall y(y \in \dot{x} \rightarrow y \cup \{y\} \in \dot{x})) \ \& \  E(x,IN,\phi) \ \ [x \ denoted \ by \ \omega]$
\newline The Sum Set makes sure that all integers are explicated.
\item Separation \newline
$\forall t_1,...,t_k \ \forall x \exists \dot{y} (\forall  z(z \in \dot{y} \iff (z \in x\   \& \   
A_n(z, t_1, ...,t_k)))) \rightarrow E(y,SEP,x,n)      \ $
\newline For infinite subclasses, if they can't Separate, they will be filtered out and will not become sets.
\item Replacement \newline
$\forall t_1,...,t_k \ \forall x \ \exists \  !\  y \ (A_n(x,y; t_1,...,t_k)) \rightarrow \forall u \exists \dot{v} \exists \dot{f} \& 
\newline \langle u,\dot{v}\rangle \in \dot{f} \  \forall r \  (r \in \dot{v} \iff \exists s \  (s \in u \ \& A_n(s,r; t_1,...,t_k) )\rightarrow 
\newline \langle s,r \rangle \in \dot{f}) \rightarrow (E(v,REP,x,y,u,n) \& E(f,REP,x,y,u,n))$ \newline
Replacement implies Separation. The $\langle u,v \rangle$s are called set-explicating functions (hereafter SEF). They do not need to be 1-1.
\item Power Set\newline
$\forall x \ \exists \dot{y} \ (\forall z \ (z \  \in \dot{y} \iff z \subseteq x \    ))\rightarrow E(y,EPS,x) \ \ [y\  is \ denoted \ by \ EPS(x)]$
\begin{enumerate}
\item ``Not a Set'' Axioms
\newline Because in EX sets have an ``only-if'' condition, clarity can sometimes be added by saying what is \textbf{not} a set. Note that these ``not a set axioms'' don't cause anything to be explicated.
\newline $ (x \ infinite) \rightarrow \exists \dot{z} (\dot{z} \subset x  \  \& \  \dot{z} \notin EPS(x)) $. 
\newline In fact, the intuition is that there are an uncountable number of such SSCs, excluded from explication by the fact that their members can't be given in a finite statement. Filtered out once in Separation and again here in the Power Set.
\item Fate of $P(\omega)$
\newline We remark at this point that the power class of $\omega$ is \textbf{not} a set in EX. (As SSCs, $EPS(\omega) \subset P(\omega))$.   
$EPS(\omega)$ must assume the role of the continuum(C) in EX. 
\newline It is easy, looking at the Axion of Specification, to see that $EPS(\omega)$ contains the normal model of the real numbers.  That is an infinite subset of $\omega$, controlled by a function, so that each member of the subset is specified.
%
\end {enumerate}
\item Rules for Explication \newline
$ E(x) \iff (x=\phi)  \lor (x=\omega) \lor (int(x)) \lor (E(x,UP,y,z)) \lor (E(x,SS,y)) \lor  (E(x,SEP,y,n)) \lor ((x,EPS,y)) \lor $
\newline $  (E(x,REP,y,z,u,n)) \lor (E(f,REP,x,y,u,n)) $  
\newline And, of course, $ M(x) \iff E(x)$
\end{enumerate}
\section {Axiom of Limited Well-Ordering}
There is No Explicated Well-Ordering For EPS of Any Infinite Set. 

\noindent Axiom - $ (x \ infinite) \rightarrow \forall y \ \neg(y \ Well-Orders(EPS(x))) $   \begin {enumerate}
\item Corollary - The Well-ordering Theorem is false in EX.  \newline
If AC were true in EX, $EPS(\omega)$ could be well-ordered, thus
\item Corollary - AC is false in EX. This allows us to filter out the unwanted, ``class-sized'' number of sets associated with AC being true, none of whose members, of course, can be given in a finite statement.
\item Corollary - No SEF maps $EPS(\omega)$ 1-1 into \textbf{any} ordinal, 
\newline including $\aleph_1$. Thus, 
\item Corollary - CH is false in EX, thus so is GCH, thus
\item Corollary - Since both AX and GCH are false, so is V=L. \newline See Appendix B for further 
discussion.
\item Discussion - This axion is extremely well-grounded in intuition, as anyone who has spent time seeking a well-ordering for  $ P(\omega) $ will attest. In fact, AC was added to ZF for just this reason.  As ZF developed, it was noticed that $P(\omega)$ was nothing like the power sets of finite sets.  Instead of being naturally well-ordered, it strongly resisted well-ordering at all. Three possible choices for resolving this dilemma were: a) hold your nose to its lack of intuitive grounding and add AC; b) instead of AC, add the even less intuitive CH as an axiom; or, c) accept that $P(\omega)$ and the ``smaller'' $EPS(\omega)$ can not be well-ordered, and get along with the development of the theory.  Choice a) was the one chosen by ZF; choice c) was taken by EX.
\end {enumerate}
\section{Developments, Conjectures and Sketches in EX}
\begin{enumerate}
\item  Regularity \newline
The axiom of Regularity is conjectured to be a theorem. The proof would start with the last explication sequence before regularity broke, then show that adding one more explicating axiom would not break it.
\item $ EPS(\omega) $ is Uncountable in EX \newline
This follows immediately from Corollary 3 of Section 4.  It also is shown by the usual approach of Cantor's Theorem ``set of all sets not members of the range of the assumed function from $ \omega$ onto $EPS(\omega)$''.
\item The Usual Higher Well-ordered Cardinals Exist in EX \newline
The usual proof goes through in $EPS^4(\omega)$  with the unordered pairs, ordered pairs, partial or full-well-orderings, and equivalence classes of well-orderings by their ordinal length, combining to produce the set of all countable ordinals ($\aleph_1$).
\item Two Dimensions of Cardinality 
\newline EX thus has 2 ``Dimensions'' of cardinality. The well-ordered ones from $EPS^4(\omega)$ and the others from $EPS(\omega)$. This doesn't cause any problems, if one remembers the cardinalities of the two dimensions don't overlap.
\item  Different Models for ZF-AC \newline
Appendix B shows that the Constructable Sets (V=L) are a model for ZF-AC.  In that model, both AC and GCH are theorems.

\noindent The explicated sets (EX) are also a Model for ZF-AC, quite different from V=L.  
This is essentially clear, as each Axiom in ZF, except AC, corresponds to an Axion in EX, but for explicated sets only. This model is interesting for two reasons. First, AC is false in EX; second, any other model of ZF (including ZF itself) with or without AC, must contain all the sets in EX.  
\item EX Decides several Problems which Took a lot of Thinking in ZF. \newline
Besides Regularity and a better solution for AC, we have \newline
$ \neg AC \rightarrow (\neg CH \ \& \neg GCH \ \& \neg (V=L))$. Also it's obvious in EX that there is no unreachable cardinal. All in all, EX is a cleaner, more-intuitive way to think about sets.
\end{enumerate}
\section{Appendix A - G.B. Axioms for Class Formations}
\begin{enumerate}
\item $\exists X \ \forall a(a \in X\longleftrightarrow \exists b \ \exists c(a = \langle b,c\rangle \ \& \ b \in c))$
\item $\forall X \ \forall Y \ \exists Z \forall u(u \in Z \longleftrightarrow u \in X \ \& \ u \in Y)$
\item $\forall X \ \exists Y \ \forall u(u \in Y \longleftrightarrow \ \sim u \in X)$
\item $\forall X \ \exists Y \ \forall u(u \in Y \longleftrightarrow \exists v( \ \langle v,u \rangle \ \in X))$
\item $\forall X \ \exists Y \ \forall u(u \in Y \longleftrightarrow \exists r \ \exists s(u = \langle r,s \rangle \ \& \ s \in X))$
\item $\forall X \ \exists Y \ \forall a(a \in Y \longleftrightarrow \exists b,c( \ \langle b,c \rangle = a \ \& \  \langle c,b \rangle \ \in X))$
\item $\forall X \ \exists Y \ \forall u(u \in Y \longleftrightarrow \exists a,b,c( \langle a,b,c \rangle \in \ X \ \& \ \langle b,c,a \rangle \in Y \ \& \ \langle b,c,a \rangle = \nolinebreak[4] u))$
\item $\forall X \ \exists Y \ \forall u(u \in Y \longleftrightarrow \exists a,b,c( \langle a,b,c \rangle \in \ X \ \& \ \langle a,c,b \rangle \in y \ \& \ \langle a,c,b \rangle = \nolinebreak[4] u))$ 
\end{enumerate}
\section {Appendix B - $Consis(ZF + AC + GCH)$}
In 1958 G\"{o}del proved that if ZF (without AC) is Consistent, it remains so if AC and GCH are added as Axioms. The blow by blow is given in COHEN[1], pp 85-99. Briefly, he showed that the Constructible Sets (L) are a class-sized model for ZF minus AC. In that 
Model(V=L), both GCH and AC are theorems.  He uses a proof principle called ``Trans-Finite Induction'', where each ordinal ``stage'' depends on the power class of the previous stage, or SUP of the previous stages for limit ordinals.  Note that in ZF all axioms are ``set-defining''. EX has the ``no-such-explicated-set-type'' Axiom of Limited Well-Ordering. That axiom would need to be true in any model of all EX's axioms. Though not relavent, one could attempt an ``Every Set is Constructible'' model using only the set-explicating axioms of EX.  Even that plan bogs down, because: 1) In EX, Trans-Finite Induction does not work to class-sized ordinals, at most only within a very large set; and 2) the $\omega th$ stage of the construction fails because $P(\omega)$ is not a set.
\section {References}
[1] P. J. COHEN: "Set Theory and the Continuum Hypothesis," Stanford University, 1966, W. A. Benjamin, Inc. Advanced Book Program, Reading, Massachusetts.    
\section{End of Article: Set Theory Axioms Using "Explication"}
\end{document}